%From shelah@math.huji.ac.il  Sun May 16 10:32:29 1999
%Date: Sun, 16 May 1999 10:32:27 +0300 (IDT)
%From: Saharon Shelah <shelah@math.huji.ac.il>
%To: Mijrady biyrujalyim <shlhetal@sunset.ma.huji.ac.il>,
%        CIGDEM GENCER <gencer@math.huji.ac.il>
%Subject: Plain Tex  698  !
%Message-ID:
%<Pine.SUN.3.95-heb-2.07.990516103129.21493B-100000@sunset.ma.huji.ac.il> 

%\batchmode
\magnification\magstephalf
\def\cf{{\rm cf}}
\font\fr=eusb10
\newfam\fra
\textfont\fra=\fr
\def\ccal{{\fam\fra C}}
\def\st{ such that }
\def\Pf{ Proof }
\def\seq{ sequence }
\def\cont{ continuous }

\def\sta{ stationary}
\def\stat{ stationary }
\def\shelah{ \centerline{Saharon shelah}
\centerline{Institute of Math, The Hebrew University, Jerusalem
Israel}
\centerline{Department of Math., Rutgers University, New
Brunswick NJ USA }}

\def\a{\alpha}
\def\A{\Rightarrow}

\def\b{\beta}

\font\mmm=msam10
\textfont9=\mmm
\mathchardef\Bo="3911

\def\c{\subseteq}

\def\d{\delta}

\def\e{\epsilon}

\def\g{\gamma}

\def\k{\kappa}

\def\l{\lambda}

\def\m{\mu}

\def\p{ \pi }

\def\s{\sigma}

\def\th{\theta}

\def\w{\omega}

\def\x{\chi}

\def\z{\zeta}
\def\phi{\varphi}

\def\p0{\emptyset}
\def\pa{\forall}
\def\pb{\setminus}
\def\pc{\wedge}

\def\pL{||}

\def\pN{\not=}
\def\po{\emptyset}
\def\ps{\subseteq}

\def\pu{\cup}
\def\pU{\bigcup}
\def\pv{\cap}

\def\py{\exists}

\def\mathcal#1{{\cal #1}}

\def\sC{\mathcal{C}}

\def\sE{\mathcal{E}}

\def\sH{\mathcal{H}}

\def\sM{\mathcal{M}}

\def\sN{\mathcal{N}}

\def\sP{\mathcal{P}}

\def\sU{\mathcal{U}}

\def\0{emptyset}
\def\ha{\aleph}

\def\lk{\langle}
\def\rk{\rangle}

\def\df{=^{def}}

\def\then{{\bf then\ }}
\def\Then{{\bf Then\ }}
% big parentheses
\def\({\big(}
\def\){\big)}
\def\[{\big[}
\def\]{\big]}

\def\0{\emptyset}
% The ``name'' macro put a \sim sign(big tilde) under its argument
\def\name#1{\mathpalette\putsimunder{#1}}
\def\putsimunder#1#2{\oalign{$#1#2$\crcr\hidewidth
\vbox to.2ex{\hbox{$#1\widetilde{\hphantom{#2}}$}\vss}\hidewidth}}

\def\no{\noindent}
\def\nl{\hfill\break}

\centerline{\bf On the existence of large subsets of
$[\l]^{<\k}$}
\centerline{{\bf which contain no unbounded non-stationary
subsets}
\footnote*{Publication number 698 in author's list.\hfill\break 
Partially supported by the Israel Science Foundation, founded by the Israel
Academy of Sciences and Humanities}}
\medskip

\shelah
\medskip

Here we deal with some problems posed by Matet.
The first section deals with the existence
of stationary subsets of $[\l]^{<\k}$
with no unbounded subsets which are not stationary, where,
of course, $\k$ is regular uncountable $\le \l$.
In the second section  we deal with the existence of such clubs.
The proofs are easy but the result seems to be very surprising.

Theorem 1.2 was proved some time ago by Baumgartner (see Theorem 2.3 of
[Jo88]) and is presented here for the sake of completeness.
\bigskip

\centerline{\bf Section 1: On stationary sets with no unbounded
stationary subsets}
\medskip

The following quite answers the question of Matet indicated in the title
of this section.
Theorems 1.2 and 1.4 yield existence under reasonable
cardinal-arithmetical assumptions.
Theorem 1.5 yields consistency.
Remark 1.3 recalls that under the assumption
$\l=\l^{<\k}$ there is no such set.
\medskip

\no1.0 NOTATION: $\sH(\x)$ is the family of sets $x$ such that ${\rm
TC}(x)$, the transitive closure of $x$, has cardinality $<\x$, and
$<^*_\x$ is any well ordering of $\sH(\x)$.
\medskip 

\no1.1 DEFINITION:
(1) For   $\k$  regular uncountable
and $\k \le \l$, let $\sE_{\l,\k}$, the club filter of 
$[\l]^{<\k}$, be the filter generated by the clubs
of $[\l]^{< \k}$, where a club of $[\l]^{< \k}$ (or
an $\sE_{\l,\k}$-club) is a set of the form
$\sC_{\sM} \df \{ a\in [\l]^{< \k}  : a = cl_{\sM} (a)
\pc  a \pv \k \in \k \}$ 
%is a club of $[\l]^{< \k} $,
where $\sM$ is a model with universe $\l$ and countable vocabulary.
Note that we could have required in the definition of $\sC_{\sM}$ that
$\sM$ restricted to $a$ be an elementary submodel of $\sM$, but it
does not matter as we can expand $\sM$ by Skolem functions.
Note that, as a consequence, only the functions of $\sM$ matter.\nl
(2)  We call $S \ps [\l]^{<\k}$ stationary if it's complement
$[\l]^{<\k} \pb S$ does not belongs to
$\sE_{\l,\k}$.\nl
(3) We call $S \ps [\l]^{<\k}$  unbounded if every member of
$[\l]^{<\k}$ is contained in some member of $S$,
so $\cf \( [\l]^{<\k}, \ps \)$ is exactly $Min\{ |
\sU | : \sU \ps [\l]^{<\k} \hbox{ is unbounded}\}$.\nl
(4) In part (1) we can replace $\l$ by any set $A$
which includes $\k$.
\medskip

We deal first with a special case.
\smallskip 

\no1.2 THEOREM: if  $\l=\aleph_\w, \k= (2^{\aleph_0})^+ <\l$
and $ 2^\l = \l^{\aleph_0}$
\then  there is a stationary subset of
$[\l]^{<\k}$  \st every unbounded subset
of it is stationary.
\smallbreak

\no\Pf : Let  $\lk \sM_\a:\a<2^\l \rk $ list all models
with universe $\l$ and countable vocabulary (for the sake of simplicity
fix the vocabulary, having $\aleph_0$ many $n$-place function
symbols for each $n$).
Let $\lk a_\a: \a< 2^\l \rk$ list without repetitions the countable subsets
of $\l$.

Now we define $S$ as the set of subsets $Y$
of $\l$ of cardinality
continuum \st for every  $\a <2^\l$, if $a_\a$ is 
a subset of $Y$ then $Y$ is the universe of 
an elementary submodel of $\sM_\a$.
Clearly in $\w_1$ steps we can "catch our tail", so proving
that $S$ is stationary  (in fact if $\lk N_\z:\z <\th \rk$
is an increasing \cont \seq of elementary submodels of
$(\sH(\x) , \in , <^*_\x  )$, to which
$\lk \sM_\a:\a<2^\l \rk $  and $ \lk a_\a: \a< 2^\l \rk$ belong,
each $N_\z$ is of the cardinality of the continuum, $\th$ a regular
uncountable $\le 2^{\aleph_0}$
and $\lk  N_\e:\e\le \z  \rk \in N_{\z+1}$ for $\z<\th$
(note that $ N_{1+\z}  \pv \k \in \k$ follows in this case
and for limit $\z$ in all cases), then $\pU_{\z<\th} N_\z \pv \l \in S$.
We denote by  $\sE_{\l,\k,\th}$ the filter generated by the family
$\big\{\{\pU_{\z<\th} N_\z \pv \l : \lk N_\z : \z<\th \rk$ is an
increasing sequence of elementary submodels of 
$(\sH(\x),\in,<^*_\x)$ such that  $x\in N_0\,,\,||N_\z||<\k$
and  $\lk N_\e : \e<\z \rk \in N_{\z+1} :\x>\l
\wedge x \in \sH(\x)\big\}$.
This filter is normal, and the set $S$ defined above belongs to
$\sE_{\l,\x,\th}$.
(On this filter see  [Sh 52, section 3]).

So S is stationary.
Now suppose that $S_1$ is an unbounded
subset of S which is not stationary, so it is disjoint to
some club, hence for some $\a<2^\l$, no member of $S_1$
is the universe of an elementary submodel of $\sM_\a$.
By the unboundedness assumption, for some member $Y$ of $S_1$
\ $ a_\a$ is a subset of $Y$, but, by the definition of $S$, $Y$ is
necessarily the universe of an elementary submodel of
$\sM_\a$, which contradicts the choice of $\a$.

QED

${}$

\no1.3. REMARKS:(1) See more in 2.4(2).

\no(2) Of course, it is known that  if $\l=\l^{<\k}$
then there is no subset  of $[\l]^{<\k}$ as claimed in 1.2.

[Why? let $\{a_\a:\a<\l\}$ be a listing of $[\l]^{<\k}$.
%\st $\a \notin a_\a$;
Now, if $S$ is a \stat subset of  $[\l]^{<\k}$,
for each $\a<\l$ choose a set 
$b_\a \in S$  which includes $a_\a  \pu \{\a\}$
and $b_\a \notin \{ b_\g:\g<\a\}$

Now $S_1 \df \{ b_\a:\a<\l\}$ is an unbounded subset
of  $S \ps [\l]^{<\k}$ as $a_\a \ps b_\a$
and is not \sta, as the mapping $b_\a \mapsto \a $
is a one to one choice function 
(contradicting normality).]

${}$

\no1.4 THEOREM:
Assume that $\l$ is a strong limit singular
of cofinality $\s$, $\k$ is regular and $\s < \k <\l$,
\then there is a stationary subset $S$
of $[\l]^{<\k}$
 \st every unbounded subset of it is stationary.

%(2) Assume that $\ha_0 <\k  \le

${}$

\no\Pf: Let $\lk \sM^*_\g : \g <\l \rk$ list the models
with universe an ordinal $<\l$ and countable
vocabulary, allowing partial functions.
Let $\lk \m_\e:\e<\cf(\l)\rk$ be a strictly increasing  sequence of cardinals
 with limit $\l$.
Let   $\lk \sM_\a:\a<2^\l \rk$ list all models
with universe $\l$ and  countable vocabulary.
For  $\a < 2^\l$ and $\e <  \cf(\l) $ let
$\g_{\a,\e}$  be  the  unique  ordinal
$\g$ \st $\sM_\a$ restricted to $\m_\e$  is  $\sM^*_\g$
and let  $a_\a$  be $\{
\g_{\a,\e}:  \e<\cf(\l)\}$.

Let $S$ be the set of subsets $Y$ of $\l$ of cardinality $<\k$
\st   $Y \pv \k \in \k$  and
if $\g\in Y$ then $Y$ is closed under the partial
functions of $\sM^*_\g$ .  Clearly  $S  \ps  [\l]^{<\k}$  
is closed unbounded.

It is not hard to see that if $Y \in S$ includes $a_\a$
then it is the universe  of a  
%n elementary 
submodel of $\sM_\a$  
[remember that we are using the variant Definition 1.1(1) without
elementaricity].
Hence every unbounded $S'  \ps S$ is \sta\  and we are done.

QED

${}$

\no1.5 THEOREM : Assume  that $\s < \k=\cf(\k) <
\m^{<\m}  < \l$ (the most interesting case is where
$\s= \cf(\l)$\thinspace) and that the axiom of [Sh 80] for
$\m$-complete forcing notions satisfying a strong
version of the $\m^+$-cc holds for  $\le \l^{<\k}$  dense subsets.
\Then any stationary subset of $[\l]^{<\k}$  has an unbounded subset
which is not stationary.

\no\Pf:
Let $ S \ps  [\l]^{<\k}$.
We define a forcing notion $Q=Q_S$.
A member $p$ of $Q$ is a function from a subset  of $S$ of cardinality
$<\m $ to  $\{0,1\}$  \st there is no increasing sequence
of limit length of members  of the set
$p^{-1}(1)$  with union a member of 
$p^{-1}(1)$, moreover this holds even for the union of directed systems.

Let $p \le q$ iff $q$ extends $p$ and
$x\ps y  \pc y \in p^{-1}(1) \pc x\in q^{-1}(1) \A  x \in p^{-1}(1)$.
Clearly  Q is  $\m$-complete, in fact, the union of every increasing
sequence of length $<\m$ of members of $Q$ is in $Q$.
Also, it is easy to see that $Q$ satisfies $\m^+$-chain condition.
For this proof it suffices to use $\Delta$-systems hence also the strong
version of the $\m^+$-chain condition (of [Sh:80]) holds.

Let $\name{S}_1$ be $\pU \{p^{-1}(1): p\in \name{G}_Q\}$,
it is forced to be a subset  of $[\l]^{<\k}$.
Clearly $\name{S}_1$ is not stationary as no non trivial union of
directed systems of members of it belongs to it.
Also, it is unbounded because the density claim is obviously satisfied.
Lastly we apply the axiom of [Sh:80] to this forcing notion and we are done.
(Alternatively, define a forcing notion $Q'=Q'_S$, where a member $p$
is a one to one choice function with domain a  subset  of $S$ of cardinality
$<\m $ ordered by inclusion.)

QED
\bigskip

\centerline{\bf Section 2 : On clubs containing no unbounded \stat
subsets}

${}$

Matet has further asked  whether any club of
$[\l]^{<\k}$ contains an unbounded non stationary subset.
Below we answer this quite completely,
we shall later deal with other variants and get similar results.

${}$

\no2.1. THEOREM: Suppose that $\k$ is regular
uncountable
and $\k \le \l$.\nl
If $\cf(\l) \ge  \k $  \then  every club of  $[\l]^{< \k}$ contains
an unbounded non-stationary subset.

${}$

\no\Pf: 
Let $\sM$ be a model with universe $\l$ and countable 
vocabulary, so $\sC = \sC_{\sM}
\df \{ a\in [\l]^{< \k}  : a = cl_{\sM} (a)
\pc  a \pv \k \in \k \}$ is a club of $[\l]^{< \k} $
and every club of $[\l]^{< \k} $ has this form,
so it is enough to find an unbounded
non-stationary subset of $\sC_\sM$.

We choose an increasing continuous
sequence $\lk \sM_\a : \a <\cf(\l) \rk$ of elementary
sub-models of $\sM$, each of cardinality $<\l$ and 
$\sM_\a \pv \k \in (\k +1 )$.

Let $\sU$ be the set of $ a \in \sC$ satisfying :
for some $ j<\cf(\l)$ \ $a$ is included in
$\sM_{j+2}$ and has a member in $\sM_{j+2} \pb \sM_{j+1}$ .
We shall now see that $\sU$ is unbounded  (as a subset of
$[\l]^{<\k}$): if $b\in [\l]^{<\k}$ then since $\cf(\l)\ge\k$
clearly for some $j<\cf(\l)$ we have $b\ps\sM_j$, let
$\a = Min\( \sM_{j+2} \setminus \sM_{j+1} \)$ and we can find a member
$a$ of~$\sC_{\sM_{j+2}}$ which includes $b \pu \{\a\}$.
So it is enough to prove that $\sU$ is not stationary,
%hence have a club of $a\pv \sM_{j+1} \pb \sM_j  \pN  \p0 $
that is, to find a club disjoint to it.
Such a club is the set $\sC^* $ of all $a\in \sC$ \st for every $j< \cf(\l)$
if $a\pv \sM_{j+1} \pb \sM_j  \pN  \p0 $
then $a\pv \sM_{j+2} \pb \sM_{j+1}  \pN  \p0 $.

QED

${}$

\no2.2. THEOREM: Suppose that $\k$ is regular uncountable
and $\l \ge \k$.

 Assume $\cf(\l) <\k$ and either $\l$ is strong limit
or at least for every $\m<\l$ there is a family
of $\le \l$ clubs of   $[\m]^{< \k} $ \st any other club
of $[\m]^{< \k} $ contains one of them.

\Then some club of $[\l]^{< \k} $
contains no  unbounded non stationary subset of $[\l]^{< \k} $.

\no\Pf:
 By the assumption we can find  $\lk (\a_\b, \sC_\b): \b < \l \rk$
\st : 

\no(a)  $\a_\b$ is an ordinal $\le \k+\b <\l$ but $\ge \k$

\no(b) $\sC_\b$ is a club of $[\a_\b ]^{<\k}$

\no(c) if $\a<\l$  and $\sC$ is a club of $[\a]^{<\k}$
then for some $\b<\l$ we have $ \a_\b = \a$ and 
$\sC_\b \ps \sC$.

Let $ \sC^* \df \{ a\in [\l]^{<\k} :
(\pa \b \in a) [ a \pv \a_\b \in \sC_\b]\}$.

Clearly, $\sC^*$ is a club of $[\l]^{<\k}$.
Towards contradiction, assume that $\sU$ is an unbounded
subset of of $\sC^*$ and we shall prove that it is stationary.
So let $\sN$ be a model with universe $\l$  and countable vocabulary,
and we shall find $ a \in \sU \pv \sC_{\sN} $ (see the proof of 2.1).
Recall that by the assumption of the theorem $\s \df \cf(\l)$ is $<\k$
hence we can find  an increasing  sequence
$\lk \g(\z) : \z <\s \rk $  of ordinals $<\l$ with limit $\l$.
Let $\sC^\z $  be $\{ a \in [{\g(\z)}]^{<\k} : 
cl_{\sN}(a) \pv {\g(\z)}  = a  \} $, clearly it is a club of
$[{\g(\z)}]^{<\k}$  hence for some  $\b(\z) <\l$ we have
$\a_{\b(\z)} = \g(\z)$  and $ \sC_{\b(\z)} \ps  \sC^\z$.

Now every member of $\sC^*$ which includes $\{ {\b(\z)} : \z <\s \}$
is necessarily the universe of a
%n elementary
submodel of $\sN$  (just think a minute), but $\sU$
contains a member $a$ which includes such a set as $\sU$ is
an unbounded of $[\l]^{<\k}$, and $a$ necessarily belongs to   $\sC_{\sN}$
as $\sU$ is a subset of $\sC_{\sN}$, so we are done.

QED

${}$

\no2.3. THEOREM: Assume $\aleph_0  < \m = \m^{<\m} < \l$
and $\l$ is a strong limit singular
and $\m  + \cf(\l) <\k=\cf(\k)  <\l <\chi$
and $P$ is any $\k$-c.c. forcing (e.g., adding $\x$~Cohen subsets to $\m$.)
\Then in $V^P$ we have:

\no(a) The cardinal arithmetic is the obvious one, and

\no(b) There is a club subset $S$ of $[\l]^{<\k}$
with no  unbounded non-stationary subset.
 
${}$ 

\no\Pf:
The point is that the condition in 2.2 continues to hold
as any new club of $[\th]^{<\k}$ contains an old one
when $ \k \le \th < \l$.

QED

${}$

\no2.4 THEOREM:
(1) Assume that $\ha_0 <\k = \cf(\k) \le \l$ and
$\l = \cf([\l]^{<\k} , \ps ) $.
\Then every \sta\ subset of
$[\l]^{<\k}$ contains an unbounded 
non-stationary
subset. 

\no(2) Assume that $\s<\k = \cf(\k) \le \l$
and $\cf([\l]^\s , \ps ) = gen(\sE^\s_{\l,\k}  )$ (see below).
\Then  some $\sE^\s_{\l,\k}$-club
$\sC^*$ contains no  unbounded 
non-stationary subset.

\no(3) Assume that $\ha_0 <\k = \cf(\k) \le \l$ and $\sE$ a filter
on $[\l]^{<\k}$ extending  $\sE_{\l,\k}$
and $\m = gen(\sE ) $ 
and $\lk a_\z :\z <\m \rk$ is a \seq of subsets of $\l
$ each of cardinality $<\k$ and for every $\x>\l,x \in \sH (\x)$
the set $\{ N \pv \l :  N \prec (\sH ( \x ) , \in ),
x \in N , \pL N \pL <\k, N \pv \k \in \k, 
(\pa \z<\m)[ a_\z \ps N \A \z \in N ] \}$ belongs to $\sE$.
\Then some member  $\sC^*$  of $\sE$
contains no  unbounded subset  which is $= \po\,{\rm mod}\,\sE $.

\no(4) Assume that $\ha_0 <\k = \cf(\k) \le \l$ and $\sE  \ps\sE'$ are filters
on $[\l]^{<\k}$ extending  $\sE_{\l,\k}$
and $\lk \sC_\z : \z < \m \rk$ list a subset of $\sE'$
which generates a filter extending $\sE$
and $\lk a_\z :\z <\m \rk$ is a \seq of subsets of $\l
$ each of cardinality $<\k$ and for every $\x>\l,x \in \sH (\x)$
the set $\{ N \pv \l :  N \prec (\sH ( \x ) , \in ),
x \in N , \pL N \pL <\k, N \pv \k \in \k, 
(\pa \z<\m)[ a_\z \ps N   \A N \pv \l \in \sC_\z  ] \}$ belongs to
$\sE'$. 

\no\Then  some member  $\sC^*$  of $\sE'$
contains no  unbounded subset  which is $= \po\, mod\, \sE $.

We shall prove 2.4 below.

${}$

\no2.5 DEFINITION :(1) $\sE^\s_{\l,\k}$ is the filter on $[\l]^{<\k}$
generated by the  $\sE^\s_{\l,\k}$-clubs, where 

\no(2) An  $\sE^\s_{\l,\k}$-club is, for some $\x>\l$ and
$ x \in \sH(\x)$, the set\nl
$\sC^\s_{\l,\k}[\x, x] \df \{ N \pv \l : N \prec (\sH(\x),\in,<^*_\x ),
\pL N \pL < \k , N \pv \k \in \k, 
(\pa X \in [N]^\s ) [ X \in N ], x\in N \}$\hfill\break
(Actually, it is clear that we can assume, without loss of generality,
that $\x=(\l^\k)^+$.
For this definition to be meaningful, i.e., for the filter to be non
trivial, we have to assume that\break
$(\forall\alpha<\kappa)(|\alpha|^\sigma<\kappa)$ ).

\no(3) For a filter $D$ let $gen(D)$ be the minimal cardinality
of s subset which generates $D$, this is the same as the cofinality
of $D$ under inverse inclusion.

${}$
 
We note the following easy monotonicity properties:
\smallskip

\no2.6 OBSERVATION:
(1) If for some filter $\sE$ extending $\sE_{\l,\k}$
we have "there is a member of $\sE$  (or just a member
of $\sE^+$) \st every unbounded
subset belongs to $\sE^+ $ ",
\then
there is a \sta\ subset of $[\l]^{<\k}$
\st every unbounded subset of it is \stat.

\no(2) Similarly we can replace $\sE_{\l,\k}$ by any filter
$\sE' \supseteq \sE $.

\no(3) Concerning 2.4(3),(4) the last assumption holds if, e.~g.,
$\l={\rm \cf}(\l)$, \ $S\c\{\d<\l:\cf(\d)<\k\}$ is \sta\ not reflecting
in any $\d<\l$ of cofinality $\k$.

${}$

\no2.7 REMARK: (1) So for $\s=\ha_0$ we get that
  $\sE^\s_{\l,\k} = \sE_{\l,\k}$.\nl
(2) Note that, of course, for every $\l$ we have
$2^\l \ge gen(\sE_{\l,\k})$ and  for
$\l$ strong limit of cofinality $<\k$   we
have
$gen(\sE_{\l,\k} ) = 2^\l$

${}$

\no\Pf  of 2.4:

\no(1) As in 1.3(2).

\no(2)  Let $\m  \df \cf([\l]^\s , \ps )$  and
$\lk \sC_\z=  \sC^\s_{\l,\k} [ \x_\z,  x_\z ] :\z<\m  \rk$
list a family which generates $\sE^\s_{\l,\k}$,
and let $\x^* \df \sup\{2^{\x_\z^+} : \z <\m \}$ and $\x=2^{2^\x}$.

Let $\lk a_\z : \z<\m \rk$ list with no repetition a cofinal subset of 
$[\l]^\s$ of cardinality $\m$, lastly let $y  \in  \sH(\x)$
code $\s,\k,\l, \m$ and the three sequences 
$\lk\x_\z :\z <\m \rk$, $\lk x_\z):\z <\m \rk$,
$\lk a_\z:\z <\m \rk$ and let $\sC^* \df
\sC^\s_{\l,\k} [\x , y ]$.

Clearly $\sC^*$ is a $\sE^\s_{\l,\k}$-club, so to finish
assume that $\sU$ is an unbounded subset of $\sC^*$
(that is, cofinal in $[\l]^{<\k}$) which is not \sta,
hence there is $\z<\m$ \st  $\sC_\z $ is disjoint to $\sU$.
As $\sU$ is unbounded, there is $a \in \sU$
\st  $a_\z \ps a $. But $\sU \ps \sC^*$ hence
$a \in  \sC^*$ so there is a model $N$ as in the definition of
$\sC^*$ (see [2.5(2)) which  witnesses $a\in \sC^*$ so in particular
 $a \df N \pv \l   \in \sC^*$.
As $a_\z \ps a \ps N$ and $| a_\z | <\s$ (by the choice of
$\lk a_\xi :\xi <\m \rk$ ), necessarily we have
$ a_\z \in  N$ but $\lk a_\z : \z<\m \rk \in N $
is without repetitions, hence $\z \in N$, and since
the sequences $ \lk \x_\z:\z<\m\rk , \lk x_\z:\z<\m\rk$
belong to $N$ clearly 
$\x_\z , x_\z \in N $ hence $N$ witnesses that
$N \pv \l \in \sC_\z$  so $a \in \sC_\z$
contradicting $a \in \sU$ and the choice of $\z$.

\no(3),(4) Proof similar to Part (2).

QED

${}$

The following gives somewhat more than
``Theorem 2.2 applied for the club filter still holds under weaker version
where $\{a\in [\l]^{< \k} : a\pv \k \in \k \}$ is not required to belong
to the filter''.

${}$

\no2.8.  DEFINITION:
(1) We say that $\ccal$ is a club filter $(\l,\k)$-definor
if  it is a function \st:

\no(a) Its domain is the family of models $\sM$ with universe
a subset of $\l$ which includes $\k$ and with countable vocabulary,
$\ccal_\sM $ is a family of submodels of $\sM$ of cardinality $<\k$.

\no(b) If $\sM$ is a submodel of $\sN$ then $\ccal_\sM  \ps \ccal_\sN$.

\no(c) If $\sM,\sN$ are models as above and $\sM$ is interpretable in
$\sN$ with the parameters from $A$ then $a\in \ccal_\sN  \pc A\ps  a$
implies $a\in \ccal_\sM$.

\no(d) If $\sM$ is as above, and $\sM_1$ is a submodel
of $\sM$ which includes $\k$ then $\ccal_\sM $ has a member
$a$ which is contained in $\sM_1$.

\no(2) For  $\ccal$ as above and  a subset $Y$ of $\l$ which includes
$\k$, the filter  $\sE [ Y ]=\sE[Y,\ccal]$  is the family of subsets of
$[Y]^{<\k}$
which include some $\ccal_\sM $ with $\sM$ as above with universe $Y$;
any such set is called an $\sE [Y]$-club  of
 $[Y]^{<\k}$. If $Y=\l$ we may write just $\sE$.
If $\sM$ is a model with countable vocabulary
and universe $Y$, $\k \ps Y \ps \l$, we may write
$\sE[\sM ]$  for  $\sE [Y ]$.
 
\no(3) If for every  such $Y$ we have $[Y]^{<\s} \in \sE[Y] $
then we call $\sE$  a $(\l,\k,\s)$-definor and may consider only
$[Y]^{<\s}$.

\no(4) We say $\sU \ps [Y]^{<\k}$    is $\sE$-unbounded if
$\{ a \in [Y]^{<\k} : (\py b \in \sU ) [ a \ps b ]
\} \in \sE [Y] $

\no(5) We say $\sU \ps [Y]^{<\k}$    is $\sE$-\stat if
$\sU \not= \po\,{\rm mod}\, \sE [Y] $

${}$

\no2.9. THEOREM: Suppose that $\ccal$ is a club filter
$(\l,\k)$-definor.
We shall abbreviate $\sE[\l,\ccal]$ by $\sE$.

\no(1) Assume that:  every $\sM$  as above, can be represented as a
union of an increasing
chain of submodels $\lk \sM_\z : \z < \d \rk $  
with $\d$ a limit ordinal,
\st\nl
$\{ a \in [\l]^{<\k}  : (\py \z<\d) [ a \ps  \sM_\z] \} \in \sE $.\nl
\Then every $\sE$-club of
$[\l]^{<\k}$ contains an $\sE$-unbounded non $\sE$-stationary subset.

\no(2) Assume $\cf(\l) <\k$ and $\l$ is strong limit
(or at least for every $\m<\l$ we have there is a family
of $\le \l$ \ $\sE[\m,\ccal] $-clubs of   $[\m]^{< \k} $ \st any other
$\sE[\m,\ccal]$-club
of $[\m]^{< \k} $ contains one of them,
in other words $\l\ge \cf( \sE[\m,\ccal], \supseteq)$ ),
and for any subset $Y$ of $\l$ of cardinality $\le \cf(\l)$  the
family $\{a\in  [\l]^{< \k} : Y \ps a \} $ belongs to the filter $ \sE $,
and if $\lk \m_\z: \z < \cf(\l)\rk $ is an increasing sequence of
ordinals~$>\k$ with limit $\l$ then for every $\sE$-club $\ccal$ of
$[\l]^{<\k}$ there is a sequence $\lk \sC_\z:\z<\cf(\l) \rk$ with
$\sC_\z$ being a $\sE[\m,\ccal]$-club \st
$(\pa a \in [\l]^{<\k}) (\pa \z <\cf( \l )  )[ a \pv \m_\z \in \sC_\z]$
 implies  $ a\in \sC$.

\no\Then some $\sE$-club of $[\l]^{< \k} $
contains no  $\sE$-unbounded non $\sE$-stationary subset of $[\l]^{< \k} $ .

${}$

\no\Pf: Similar to the previous ones, e.g.

\no(1) Let $\sC$ be a $\sE$-club, say $\sC = \sE_\sM$,
and let $\lk \sM_\z : \z < \d \rk $ be as above. Now
we choose\nl 
$\sU \df \{ N : \hbox{ for some }\z<\d, N \ps \sM_{\z+2}\>\hbox{and}\>
N \setminus \sM_{\z+1} \not= \po \}$.\nl
Now use clause (b).

QED

${}$

${}$

\centerline{\bf REFERENCES }

${}$

\no[Jo88] C.A. Johnson,
 Seminormal $\lambda$-generated 
ideals on  ${\sP}_\kappa\lambda$,
Journal of Symbolic Logic 53 (1988) pp. 92-102.\nl
[Sh:52] S.Shelah, 
A compactness theorem for singular cardinals, free algebras,
Whitehead problem and transversals,
Israel Journal of Mathematics, 21 (1975), pp. 319--349.\nl
[Sh:80] S.Shelah,
A weak generalization of MA to higher cardinals,
Israel Journal of Mathematics~30 (1978), pp.~297--306.
\bye